\documentclass[a4paper, reqno, 11pt]{amsart}

\usepackage[english]{babel}
\usepackage{amsmath}
\usepackage{amssymb}
\usepackage{enumerate}
\usepackage{ifthen}
\usepackage{bbm}
\usepackage{color}
\provideboolean{shownotes} 
\setboolean{shownotes}{true}
\usepackage{hyperref}

\newcommand{\margnote}[1]{
\ifthenelse{\boolean{shownotes}}%
{\marginpar{\raggedright\tiny\texttt{#1}}}%
{}%
}

\newcommand{\hole}[1]{
\ifthenelse{\boolean{shownotes}}%
{\begin{center} \fbox{ \rule {.25cm}{0cm}
\rule[-.1cm]{0cm}{.4cm} \parbox{.85\textwidth}{\begin{center}
\texttt{#1}\end{center}} \rule {.25cm}{0cm}}\end{center}}
{}
}
\newtheorem{theorem}{Theorem}[section]

\newtheorem{lemma}[theorem]{Lemma}

\newtheorem{definition}[theorem]{Definition}

\theoremstyle{remark}

\newcommand{\R}{\mathbb{R}}

\newcommand{\pt}{\partial_t}
\newcommand{\uv}{u^{\nu}}

\newcommand{\om}{\omega}
\newcommand{\omv}{\omega^{\nu}}
\newcommand{\dive}{\mathop{\mathrm {div}}}

\newcommand{\phiv}{\phi^{\nu}}
\numberwithin{equation}{section}



\begin{document}

\title[Renormalized Solutions 2D Euler Equations]{Renormalized Solutions of the 2D Euler Equations}

\author[G. Crippa]{Gianluca Crippa}
\address[G. Crippa]{
Departement Mathematik und Informatik, Universit\"at
  Basel\\ Rheinsprung 21 \\CH-4051 Basel \\Switzerland}
\email[]{\href{gianluca.crippa@unibas.ch}{gianluca.crippa@unibas.ch}}

\author[S. Spirito]{Stefano Spirito}
\address[S. Spirito]{
GSSI - Gran Sasso Science Institute\\ Viale Francesco Crispi 7\\ 67100 L'Aquila\\ Italy}
\email[]{\href{stefano.spirito@gssi.infn.it}{stefano.spirito@gssi.infn.it}}

\begin{abstract}
In this paper we prove that solutions of the 2D Euler equations in vorticity formulation obtained via vanishing viscosity approximation are renormalized.
\end{abstract}

\maketitle

\section{Introduction}
Let us consider the Euler equations for an incompressible fluid in vorticity formulation in $\R^2\times(0,T)$: 
\begin{equation}\label{eq:E}
\begin{aligned}
\partial_t\omega+u\cdot\nabla\omega&=0\\
u&=K\star \omega,
\end{aligned}
\end{equation}
where $\omega\in\R$ is the vorticity, $u\in\R^2$ is the velocity field and the second equation of \eqref{eq:E} is the Biot-Savart law, which reads as follows: 
\begin{equation}\label{eq:BS}
u(x,t)=\frac{1}{2\pi}\int_{\R^2}\frac{(x-y)^{\perp}}{|x-y|}\omega(y,t)\,dy.
\end{equation}
The system \eqref{eq:E} is as usual coupled with the initial datum 
\begin{equation}\label{eq:IDE}
\omega(x,0)=\omega_0\qquad\textrm{on}\qquad\R^2\times\{t=0\}.
\end{equation}
A classical problem in fluid dynamics is the approximation in the limit of vanishing viscosity ($\nu\to 0$) of solutions of \eqref{eq:E} by solutions of the Navier-Stokes equations, which in vorticity formulation read as follows:
\begin{equation}\label{eq:NS}
\begin{aligned}
\partial_t\omv+\uv\cdot\nabla\omv&=\nu\Delta\omv\\
\uv&=K\star\omv\\
\omv(x,0)&=\omv_0,
\end{aligned}
\end{equation}
where $\omv_0$ is a suitable approximation of $\om_0$.

In~\cite{LLM}, Proposition 1, the authors noted that {\em any} weak solution $\omega$ of \eqref{eq:E} in $L^{\infty}(0,T;L^{p}(\Omega))$ with $p\geq 2$ is renormalized in the sense of DiPerna and Lions \cite{DPL}, and in the case when $1<p<2$ all weak solutions obtained as a limit of {\em exact} smooth solutions of \eqref{eq:E} are renormalized. As pointed out in \cite{LLM} these results are a direct consequence of the results contained in Theorem II.3 of \cite{DPL}, see also Section \ref{sec:Pre} Theorem \ref{teo:T} in the following. Then, in \cite{LLM} it has been posed the question  whether solutions of \eqref{eq:E}-\eqref{eq:IDE} with $\omega$ in $L^{\infty}(0,T;L^{p}(\R^2))$ obtained in the {\em vanishing viscosity approximation} are renormalized when $1<p<2$. In this note we give a positive answer to the above question. Our main result is as follows, see Section \ref{sec:Pre} for the notations and the main definitions.
\begin{theorem}\label{teo:main}
Let $1<p<2$. Let $\om_0\in L^{p}_{c}(\R^2)$. Let $\{\omv_0\}_{\nu}\subset L_{c}^{p}(\R^2)$ be a sequence of smooth functions converging strongly to  $\om_0$ in $L^{p}(\R^2)$ and  $(\uv, \omv)$ be the unique smooth solution of \eqref{eq:NS}. Then, there exist $$(u,\om)\in L^{\infty}(0,T;W_{loc}^{1,p}(\R^2))\times L^{\infty}(0,T;L^{p}(\R^2))$$ such that 
\begin{equation*}
\begin{aligned}
&\uv\rightarrow u\textrm{ strongly in }L^{p}(0,T;L^{p}_{loc}(\R^2))\\
&\omv\stackrel{*}{\rightharpoonup} \om\textrm{ weakly* in }L^{\infty}(0,T;L^{p}(\R^2)). 
\end{aligned}
\end{equation*}
Moreover, $u=K\star\omega$ for a.e. $(x,t)\in \R^2\times(0,T)$ and the vorticity equation is satisfied in the renormalized sense.
\end{theorem}
This problem has been raised in the paper \cite{LLM} in connection with the enstrophy defect of weak solutions of \eqref{eq:E}. However, the interest in proving this result relies also in the fact that the vanishing viscosity is a very physical approximation of the Euler equations and proving that the vorticity satisfies \eqref{eq:E}-\eqref{eq:IDE} in a renormalized sense would imply {\em a posteriori} that solutions obtained via vanishing viscosity are transported in a weak sense by the characteristic of the vector field $u$, as it holds in the smooth setting, see \cite{A,DPL}. Moreover, the notion of renormalized solutions is also useful when the regularity available is not enough to state the equations in distributional sense. This is already the case for the Euler equations, indeed when the initial vorticity is in $L^{p}$ with $1<p<\frac{4}{3}$ the a priori estimates available for $u$ and $\omega$ are not enough to guarantee that the non linear term $u\,\omega$ is in $L^{1}_{t,x}$.\par Although a positive answer to the above question is highly expected due to the physical meaning of the vanishing viscosity approximation, the rigorous proof is not straightforward. Indeed,  when we assume that $\omv_0$ is uniformly bounded only in $L^{p}$ with $1<p<2$, only a uniform bound in $L^{p}$ for $\omega^{\nu}$ can be expected, and thus by classical result in singular integrals, only a uniform bound in $W^{1,p}_{loc}$ for $u^{\nu}$. In order to prove that \eqref{eq:E} is satisfied in a renormalized sense we need to prove that the following equation is satisfied in a weak sense: 
\begin{equation*}
\partial_t\beta(\omega)+(u\cdot\nabla)\beta(\omega)=0.
\end{equation*}
Then, some strong convergence for the sequence $\omega^{\nu}$ would be needed. However, the strong convergence of $\omv$ is not available in a straightforward way in our setting. Indeed, the strong convergence is usually achieved either as a {\em consequence} of the renormalization, see \cite{DPL}, or  by arguing in the Lagrangian setting, see \cite{CDL} (see also \cite{BBC} for an application of Lagrangian techniques  to the system \eqref{eq:E} in the endpoint case $p=1$ to get strong convergence of the vorticities). However, these arguments seem not to work here. We will avoid this problem by using in a deeper way the theory of renormalized solution, in particular the duality technique introduced in \cite{DPL}. 

Finally, we want to point out that the hypothesis $\omega_0\in L^{p}_{loc}(\R^2)$ in Theorem \ref{teo:main} can be relaxed to $\omega_0\in L^{1}(\R^2)\cap L^{p}(\R^2)$. In this case the proof holds with minor changes. 
\section{Preliminaries}\label{sec:Pre}
Given $\Omega\subset\R^2\times(0,T)$, the space of compactly supported smooth functions on $\Omega$ will be denoted by $\mathcal{D}(\Omega)$. 
We will denote with $L^{p}(\R^2)$ the standard Lebesgue spaces and with $\|\cdot\|_p$ their norm. Moreover, $L^p_{c}(\R^2)$ denotes the space of $L^p$ functions with compact support. The Sobolev space of $L^{p}$ functions with distributional derivatives in $L^{p}$ is denoted by $W^{1,p}(\R^2)$.  The spaces $L^{p}_{loc}(\R^2)$ and $W^{1,p}_{loc}(\R^2)$ denote the space of functions which are locally in $L^{p}$ and $W^{1,p}$, respectively. 
\subsection{Renormalized solutions of transport equations.}
In this section we recall the main notions and results about renormalized solutions of transport equations. We want to point out that the result will not be stated in the full generality but will be adapted to our setting. For a complete overview of the theory of renormalized solutions see for instance \cite{AC,AC1} and reference therein.

 Let us consider the Cauchy problem for the transport equation: 
\begin{equation}\label{eq:T}
\begin{aligned}
\pt w+b\cdot\nabla w&=0\\
w(x,0)&={w}_0(x)
\end{aligned}
\end{equation}
with $b$ divergence-free. The definition of renormalized solutions of \eqref{eq:T} is the following:
\begin{definition}\label{def:ren}
A measurable function $w$ is a renormalized solution of \eqref{eq:T} if for any $\beta\in C^1(\R)\cap L^{\infty}(\R)$ the following Cauchy problem is satisfied in the sense of distributions:
\begin{equation}\label{eq:RT}
\begin{aligned}
\partial_t\beta(w)+b\cdot\nabla\beta(w)&=0\\
\beta(\omega)(x,0)&=\beta(w_0).
\end{aligned}
\end{equation}
\end{definition}
In \cite{DPL} DiPerna and Lions proved, among other deep results, the following theorem.
In the statement $p$ and $q$ are in $(1,\infty)$ and are unrelated unless explicitly specified. 
\begin{theorem}\label{teo:T}
Let $b$ be a vector field in $\R^2$ such that 
\begin{equation}\label{eq:assb}
\begin{aligned}
&b(t,x)\in L^{1}(0,T;W^{1,p}_{loc}(\R^2)), \qquad \dive b = 0, \\
&\frac{|b(t,x)|}{1+|x|}\in L^{1}(0,T;L^{1}(\R^2))+L^{1}(0,T;L^{\infty}(\R^2)).
\end{aligned}
\end{equation}
\begin{enumerate}
\item Let $w_0\in L^{q}(\R^2)$. Then, there exists a unique renormalized solution in the sense of Definition \ref{def:ren} of the Cauchy problem \eqref{eq:T}, and it belongs to $C([0,T];L^{q}(\R^2))$.
\item If $\frac{1}{p}+\frac{1}{q}=1$ every distributional solution in $L^{\infty}(0,T;L^{q}(\R^2))$ of \eqref{eq:T}
 is a renormalized solutions and then there exists a unique distributional solution of \eqref{eq:T} in such space. 
 \end{enumerate}
\end{theorem}
A crucial tool in the proof of the Theorem \ref{teo:main} is a duality formula for the transport equations associated to the vorticity. In particular we shall use the following theorem again proved in \cite{DPL}. 
\begin{theorem}\label{teo:dual}
Let $b$ satisfy assumptions \eqref{eq:assb} and $w\in L^{\infty}(0,T;L^{p}(\R^2))$ be a renormalized solution of \eqref{eq:T}. Let 
$\phi\in L^{\infty}(0,T;L^{q}(\R^2))$ with $\frac{1}{p}+\frac{1}{q}=1$ be a renormalized solution of the following backward transport problem
\begin{equation}\label{eq:bwt}
\begin{aligned}
-\pt\phi-\dive(b\phi)&=\chi\\
\phi(T,x)=\phi_T(x)
\end{aligned}
\end{equation}
with $\chi\in L^{1}(0,T;L^{q}(\R^2))$ and $\phi_T\in L^{q}(\R^2)$. Then, the following formula holds
\begin{equation*}
\iint\chi w\,dxdt=\int\phi(x,0)w_0(x)\,dx-\int\phi_T(x)w(x,T)\,dx,
\end{equation*}
where the right-hand side makes sense because of Theorem \ref{teo:T}.
\end{theorem} 

\subsection{Vanishing viscosity limit}
In this section we recall the main results concerning the vanishing viscosity approximation of \eqref{eq:E}-\eqref{eq:IDE}. We want to point out that this is now a well-understood theory at least in the case there are no boundaries. The following results can be found for instance in the  monograph \cite{MB}.

 First, it well-know that \eqref{eq:NS} are globally well-posed in two dimensions.
\begin{theorem}\label{teo:euns}
For any smooth initial datum $\omv_0$ there exists a unique smooth solution of the Cauchy problem \eqref{eq:NS}.
\end{theorem} 
Moreover, the following uniformly with respect to $\nu$ a priori estimate holds:
\begin{lemma}\label{lem:lpest}
Let $(\uv,\omv)$ be the unique smooth solution of \eqref{eq:NS}, then for any $1\leq p\leq \infty$ 
\begin{equation}\label{eq:lpest}
\sup_{t\in(0,T)}\int|\omv(x,t)|^p\,dx\leq \int|\omv_0|^p\,dx.
\end{equation}
\end{lemma}
By using the Biot-Savart law, the Calder\'on-Zygmund theorem for singular integrals and \eqref{eq:lpest} we get the following a priori estimate  for the velocity field $\uv$
uniform with respect to $\nu$. 
\begin{lemma}\label{lem:sobest}
Let $(\uv,\omv)$ be the unique smooth solution of \eqref{eq:NS}. Then for any fixed $R>0$ and $1<p<\infty$
\begin{equation}\label{eq:sobest}
\sup_{t\in(0,T)}\int_{B_R(0)}|\uv(x,t)|^p\,dx+\int|\nabla\uv|^p\,dx\leq C(R),
\end{equation}
where $C(R)$ is depending only on the the norm of the initial datum $\omv_0$ in $L^{1}(\R^2)$ and in $L^{p}(\R^2)$.
\end{lemma}
\section{Proof of the Theorem \ref{teo:main}}\label{sec:proof}
In this section we give the proof of Theorem \ref{teo:main}. We divide the proof in several steps.\\
\\
\emph{Step 1. Existence of the limit.}\\
\\
Since $\{\omv_0\}_{\nu}$ is uniformly bounded in $L^{p}_{c}(\R^2)$ and strongly convergent to $\om_0$ in $L^{p}(\R^2)$ we have that  
\begin{equation}\label{eq:convini}
\omv_0\rightarrow\om_0\quad\textrm{strongly in }L^{r}(\R^2)\textrm{ for any $r\in[1,p]$}.
\end{equation}
By using Theorem \ref{teo:euns}, Lemma \ref{lem:lpest} and Lemma \ref{lem:sobest} we get that for any $\nu$ there exists a unique smooth solution $\omv$ of \eqref{eq:NS} such that for any $r\in [1,p]$ uniformly with respect to $\nu$ we have 
\begin{equation}\label{eq:est1}
\begin{aligned}
&\{\omv\}_{\nu}\subset L^{\infty}(0,T;L^{r}(\R^2))\\
&\{\uv\}_{\nu}\subset L^{\infty}(0,T;W_{loc}^{1,p}(\R^2)).
\end{aligned}
\end{equation}
Then, by using \eqref{eq:est1} and the fact that $\omv$ solves \eqref{eq:NS} a standard compactness argument implies that there exist $u\in L^{\infty}(0,T;W^{1,p}_{loc}(\R^2))$ such that up to subsequence not relabeled 
\begin{equation}\label{eq:convu}
\uv\rightarrow u\textrm{  strongly in  }L^{p}(0,T;L^{p}_{loc}(\R^2))
\end{equation}
and $\om\in L^{\infty}(0,T;L^{1}(\R^2)\cap L^{p}(\R^2))$ such that 
\begin{equation}\label{eq:convom}
\omv \stackrel{*}{\rightharpoonup} \om \textrm{  weakly$^*$ in  }L^{\infty}(0,T;L^{p}(\R^2)).
\end{equation}
Let $\eta\in\mathcal{D}(\R^2\times(0,T))$, then 
\begin{equation*}
\begin{aligned}
\iint u\eta\,dxdt&=\lim_{\nu\rightarrow 0}\iint\uv\eta\,dxdt\\
                            &=\lim_{\nu\rightarrow 0}\iint (K\star\omv)\eta\,dxdt\\
                            &=-\lim_{\nu\rightarrow 0}\iint\omv (K\star\eta)\,dxdt.
\end{aligned}
\end{equation*}
Because of \eqref{eq:est1} and since $\eta$ is bounded and compactly supported we have that $K\star\eta\in L^{q}(\R^2\times(0,T))$ for any $q>2$. Then, by choosing $q=\frac{p}{p-1}$ we get
\begin{equation*}
\begin{aligned}
\iint u\eta\,dxdt&=-\lim_{\nu\rightarrow 0}\iint\omv (K\star\eta)\,dxdt\\
                            &=-\iint\om(K\star\eta)\,dxdt\\
                            &=\iint(K\star\om)\eta\,dxdt,
\end{aligned}
\end{equation*}
where in the second line we have used that $\omv$ is weakly convergent in $L^{p}(\R^2\times(0,T))$. Then, by varying $\eta\in\mathcal{D}(\R^2\times(0,T))$ we have that 
\begin{equation}\label{eqE2}
u(x,t)=(K\star\omega)(x,t)\qquad\textrm{a.e. in }\R^2\times(0,T). 
\end{equation}
Now we prove that $u$ satisfies the growth condition \eqref{eq:assb}.
We decompose $u$ in the following way:
\begin{equation*}
u=K\star \omega=u_1+u_2=K_1\star\omega+K_2\star\omega
\end{equation*}
where $K_1(x)=K(x)\mathbbm{1}_{|x|\leq 1}(x)$ and $K_2(x)=K(x)\mathbbm{1}_{|x|>1}(x)$.
Then, by direct computations we have that $K_1\in L^{1}(\R^2)$ and $K_2\in L^{\infty}(\R^2)$. By using Young inequality 
and the fact that the $\om$ is bounded in $L^{1}(\R^2)$ we get \eqref{eq:assb}. \\
\\
\emph{Step 2. A dual problem}\\
\\
Let us introduce the following linear backward transport-diffusion problem 
\begin{equation}\label{eq:bns}
\begin{aligned}
-\pt\phiv-\nu\Delta\phiv-\dive(\phiv\uv)&=\chi\\
\phiv(x,T)&=0,
\end{aligned}
\end{equation}
where $\chi\in \mathcal{D}(\R^2\times(0,T))$ and $\{\uv\}_\nu$ is the subsequence chosen in \eqref{eq:convu}. By standard energy estimates it is easy to prove that for any fixed $\nu>0$ there exists a unique global smooth solution of \eqref{eq:bns}. In the remaining of this step we prove that the sequence $\{\phiv\}_\nu$ converges to the unique distributional solution $\phi\in L^{\infty}(0,T;L^{q}(\R^2))$ of the backward transport equation \eqref{eq:bwt}.

Indeed, by multiplying \eqref{eq:bns} by $2\phiv$  we get 
\begin{equation*}
-\frac{d}{dt}\|\phiv\|_2^2+2\nu\|\nabla\phiv\|_2^2=2\int\chi\phiv\,dx.
\end{equation*}
Note that the integration by parts can be justified with a suitable truncation argument as in \cite{DPL}.
Then, by using Cauchy-Schwarz inequality, the fact that $\chi$ is smooth and taking into account that the time is reversed we have by Gronwall Lemma that uniformly with respect to $\nu$
\begin{equation}\label{eq:convback}
\begin{aligned}
&\phiv\in L^{\infty}(0,T;L^{2}(\R^2))\\
&\sqrt{\nu}\nabla\phiv\in L^{2}(0,T;L^{2}(\R^2)).
\end{aligned}
\end{equation}
Moreover, by multiplying \eqref{eq:bns} by $q\phiv|\phiv|^{q-2}$ for any $q\in[2,\infty)$ we get 
\begin{equation*}
-\frac{d}{dt}\int|\phiv|^q\,dx+\nu q(q-1)\int|\nabla\phiv|^2|\phiv|^{q-2}\,dx=q\int\chi\phiv|\phiv|^{q-2}\,dx.
\end{equation*}
By using H\"older inequality, Gronwall Lemma and the fact that $\chi\in \mathcal{D}(\R^2\times(0,T))$ we get that 
\begin{equation}\label{eq:convbackp}
\phiv\in L^{\infty}(0,T;L^{q}(\R^2))
\end{equation}
uniformly with respect to $\nu$ for any $2\leq q<\infty$. Finally, by using \eqref{eq:convbackp} and the fact that $\phiv$ satisfies \eqref{eq:bns} we can improve the convergence in time, namely we have up to subsequence
\begin{equation}\label{eq:phi}
\phiv\rightarrow\phi\textrm{ in }C([0,T];L^{q}_{weak}(\R^2))
\end{equation}
for any $q<\infty$. Let $\psi\in \mathcal{D}(\R^2\times(0,T])$, by multiplying \eqref{eq:bns} by $\psi$ we have that 
\begin{equation*}
\iint\phiv\psi_t+\nu\nabla\phiv\nabla\psi+(\uv\cdot\nabla\psi)\phiv-\chi\psi\,dxdt=0.
\end{equation*}
By using \eqref{eq:convu}, \eqref{eq:convback} and \eqref{eq:phi} we get that the limiting function $\phi$ satisfies  
\begin{equation*}
\iint\phi\psi_t+(u\cdot\nabla\psi)\phi-\chi\psi\,dxdt=0.
\end{equation*}
Since by \eqref{eq:convbackp} $\phiv$ are uniformly bounded also in $L^{\infty}(0,T;L^{q}(\R^2))$ for any $q>2$ we can choose $q=\frac{p}{p-1}$ and by using Theorem \ref{teo:T} we get that $\phi$ is the unique renormalized solution of \eqref{eq:bwt}. Then, by uniqueness of the limit $\phiv$ converges to $\phi$ along the whole subsequence chosen in \eqref{eq:convu}.\\
\\
{\em Step 3. A duality formula.}\\
\\
Let us multiply \eqref{eq:NS} by $\phiv$, after integrating by parts we get 
\begin{equation*}
\iint-\pt\phiv\omv-\nu\Delta\phiv\omv-\dive(\phiv\uv)\omv\,dxdt=\int\phiv(x,0)\omv_0(x)\,dx
\end{equation*}
and by using \eqref{eq:bns} we get 
\begin{equation}\label{eq:bns1}
\iint\chi\omv\,dxdt=\int\phiv(x,0)\omv_0(x)\,dx.
\end{equation}
By using \eqref{eq:convini}, \eqref{eq:convom} and \eqref{eq:phi} we can pass to the limit in \eqref{eq:bns1} and we get 
\begin{equation}\label{eq:be1}
\iint\chi\om\,dxdt=\int\phi(x,0)\om_0(x)\,dx,
\end{equation}                                                                                                                                                                                                                                                                                                                                                                                                                                                                                                                                                                                                                                                                                                                                                                                                                                                                                                                                                                                                              
where $\phi$ is the unique renormalized solutions of \eqref{eq:bwt}.\\
\\
{\em Step 4. Renormalization.}\\
\\
By using Theorem \ref{teo:T} there exists a unique renormalized solution $\bar{w}\in L^{\infty}(0,T;L^{q}(\R^2))$ of the transport equation with vector field $u$, the limit obtained in $\eqref{eq:convu}$, and initial data $\om_0$, namely 
\begin{equation*}
\begin{aligned}
\partial_t\bar{w}+u\cdot\nabla\bar{w}&=0\\
\bar{w}(x,0)&=\om_0.
\end{aligned}
\end{equation*}
By Theorem \ref{teo:dual} we get that $\bar{w}$ satisfies 
\begin{equation}\label{eq:ed1}
\iint\chi\bar{w}\,dxdt=\int\phi(x,0)\bar{w}_0(x)\,dx
\end{equation}
where $\phi$ is the unique renormalized solution of \eqref{eq:bwt}.
By taking the difference of \eqref{eq:ed1} and \eqref{eq:be1} we get that 
\begin{equation}\label{eq:final}
\iint(\omega-\bar{w})\chi\,dxdt=0.
\end{equation}
Note that it is crucial the fact that the backward transport equation \eqref{eq:bwt} has a unique distributional solution in the class $L^{\infty}(0,T;L^{q}(\R^2))$. 
Then, by varying $\chi$ over all smooth compactly supported functions we get that $\omega-\bar{w}=0$ for a.e. $(x,t)\in\R^2\times(0,T)$ and then $\omega$ is renormalized because it agrees almost everywhere with the renormalized solution $\bar{w}$.\\

\noindent {\bf Acknowledgments.} This research has been partially supported by the SNSF grants 140232 and 156112. This work has been started while the second author was a PostDoc at the Departement Mathematik und Informatik of the Universit\"at Basel. He would like to thank the department for the hospitality and the support.

\end{document}